\documentclass[12pt]{article}
\usepackage{amsthm}
\newtheorem{Th}{Theorem}
\newtheorem{Le}{Lemma}

\usepackage{amsthm,amsmath,amssymb}

\newtheorem{Cl}{Claim}
\newtheorem{cor}{Corollary}

\DeclareMathOperator{\Fix}{\displaystyle{g}}

\def\cal#1{\mathcal{#1}}

\usepackage[all]{xy}

\usepackage[margin=2cm]{geometry}

\begin{document}
\thispagestyle{empty}
\bigskip\bigskip

\title{The solution of the complete nontrivial cycle intersection problem for permutations} \date{}
\author{Vladimir Blinovsky*, Llohann D Speran\c{c}a**}
\date{\small
 Universidade Federal de S\~ao Paulo (UNIFESP). \\
 Campus S\~ao Jose dos Campos. Instituto de Ci\^encia e Tecnologia (ICT),
 Brazil,\\
 {*}vladimir@blinovsky@unifesp.br, vblinovs@yandex.ru\\
 {**}speranca@unifesp.br, speranca@gmail.com\\ 
 \bigskip
}

\maketitle\bigskip
\newpage

\begin{center}
{\bf Abstract}

 In this paper we present a solution to the complete $t$-cycle intersection problem for systems of permutations of a finite set. 
  \end{center}

\section{Introduction and notation}


Let ${[n]\choose k }$ denote the collection of all $k$-subsets of $[n]=\{1,...,n\}$. Erd\"os--Ko--Rado \cite{erdos-ko-rado} began the investigation of the maximum size of families $\mathcal A\subseteq {[n]\choose k}$ which is $t$-intersecting, i.e., $|A\cap B|\geq t$ for any $A,B\in\mathcal A$. Later on, Deza--Frankl \cite{deza-frankl} first considered the analogous problem for permutations, along the context of coding theory. Let $\Gamma(n)$ denote the set of permutations in $[n]$. Here we follow Ku--Renshaw \cite{1} and consider  families $\mathcal A\subseteq \Gamma(n)$ which are  $t$-\textit{cycle-intersecting}, that is, for every pair $A,B\in\mathcal A$, their cycle decompositions have $t$ cycles in common. $\mathcal A$ is called a \textit{non-trivial $t$-cycle intersecting family} if there are no $t$ cycles that are common to every permutation in $\mathcal A$.
 We compute the maximal non-trivial $t$-cycle intersecting families, for every $t$.


We denote by $[r,s]$  the set of integers between $r$ and $s$, and $[n]=\{1,...,n\}$. Let  $\Gamma (n)$ be the set of permutations  of  $[n]$.  
We write $\Omega (n,t)$ for the family consisting of all systems of $t$-cycle intersecting permutations of $[n]$, and by $\tilde{\Omega} (n,t)$ the family of { systems of  pairwise nontrivial $t$-cycle intersecting permutations of $[n]$}. 

We say that $i$ is fixed by $p\in\Gamma (n)$ if $p(i)=i$ and denote 
\[f(p)=\{i\in[n]: p(i)=i \}. \]
For convenience, we often omit the word `cycle'
and refer to \textit{$t$-cycle intersecting families} of permutations as $t$-intersecting families. The aim of this paper is to determine,
\begin{eqnarray*}
&& \tilde{M}(n,t)=\max\{ |{\cal A}|: {\cal A}\in\tilde{\Omega} (n,t)\} .
\end{eqnarray*}

Before we state our main result, we present some previous results and more definitions. 

The analogous quantity $M (n,t) =\max\{ |{\cal A}|: {\cal A}\in \Omega (n,t)\}$ was studied in the literature: the case $t=1$ was settled in~\cite{4},\cite{5},
$$
M(n,1)=(n-1)!~ ,
$$
and when $n$ is large, $n>n_0 (t)$, in \cite{1}:
$$
M(n,t)=(n-t)!~.
$$

The general case was settled in \cite{3}.
 \begin{Th}[Blinovsky, \cite{3}] 
 \label{th333}
 \begin{eqnarray*}
 && M (n,t)  =\max_{r\in [0,\lfloor (n-t)/2\rfloor ]}|\{ p\in\Omega(n): |[t+2r]\cap f(p)|\geq t+r\} |.
 \end{eqnarray*}
 \end{Th}
 
 Observe that $M(n,t)$ is realized by a non-trivial $t$-cycle intersecting set whenever the maximum above ir realized for $r>0$. Moreover, as it follows from the proof of Theorem \ref{th333}, the maximum $r$ is the greatest $r$ such that $\ell=t+2r\leq n$ and
\begin{equation}
\label{eerr}
\frac{\ell -t}{2(\ell -1)}\gamma (\ell ,n )\leq 1,
\end{equation}
where 
$$
\gamma (\ell ,n)=\frac{\sum_{i=0}^{n-\ell +1} \xi\left( n-\frac{\ell +t}{2}+1-i\right){n-\ell +1\choose i}}{\sum_{i=0}^{n-\ell}\xi\left( n-\frac{\ell +t}{2}-i\right){n-\ell\choose i}}.
$$
and
$$
\xi(n)=n!\sum_{i=0}^n \frac{(-1)^i}{i!}.
$$
is the number of permutations of $[n]$ which do not have singletons.
We restate this fact as follows:
 \begin{Th}[Blinovsky, \cite{3}] 
\label{th00}
Suppose there is  $r>0$ such that $\ell =t+2r$ satisfies
\begin{equation*}
\frac{\ell -t}{2(\ell -1)}\gamma (\ell ,n )\leq 1.
\end{equation*}
Then $\tilde M(n,t)=M(n,t)$ and
\begin{equation}
\label{eel}
M(n,t)=\sum_{i=t+r}^{t+2r}{t+2r\choose i}\sum_{j=0}^{n-t-2r}{n-t-2r\choose j}\xi (n-i-j).
\end{equation}
\end{Th}

%
%

As usual, set
\[2^{[n]}=\{A\subset [n] \},\qquad\qquad {[n]\choose k}=\{A\in 2^{[n]}:|A|=k\}.\] 
Denote the set-theoretic families of $t$-intersecting systems as:
\begin{align*}
I(n,t)&=\{\mathcal A\subset 2^{[n]}: |A_1\cap A_2|\geq t,~\forall A_1,A_2\in\mathcal A \},\\
I(n,k,t)&= \{\mathcal A\in I(n,t): |A|=k,~\forall  A\in\mathcal A \},\\
\tilde I(n,t,k)&=\{\mathcal A\in I(n,k,t):|\cap_{A\in\mathcal A}A|<t \}
\end{align*}

The quantity 
$$
\tilde{M}(n,k,t)=\max_{{\cal A}\in\tilde{I}(n,k,t)}|{\cal A}|
$$
was studied by Hilton--Milner and Frankl and completely determined by Ahlswede--Khachatrian later on:
\begin{Th}[Hilton--Milner, \cite{8}]
For $n>2k$ we have,
$$
\tilde{M}(n,k,1)={n-1\choose k-1}-{n-k-1\choose k-1}+1 .
$$\end{Th}

\begin{Th}[Frankl, \cite{9}]
For sufficiently large $n>n_0 (k,t)$, 
\begin{itemize}
\item if $t+1\leq k\leq 2t+1,$ then $\tilde{M}(n,k,t)=|\nu_1 (n,k,t)|,$ where,
$$
\nu_1 (n,k,t)=\left\{ V\in{[n]\choose k}:\ |[t+2]\bigcap V|\geq t+1\right\},
$$
\item if $k>2t+1$, then $\tilde{M}(n,k,t)=|\nu_2 (n,k,t)|,$ where,
\begin{multline*}
\nu_2 (n,k,t)=\left\{ V\in{[n]\choose k}~:~ [t]\subseteq V,\ V\cap [t+1,k+1]\neq\emptyset\right\}
\bigcup \left\{ [k+1]\setminus\{ i\}:\ i\in [t]\right\}.
\end{multline*}
\end{itemize}
\end{Th}

\begin{Th}[Ahlswede--Khachatrian, \cite{7}]
\label{th66}
\mbox{}

\begin{itemize}
\item
If $2k-t<n\leq (t+1)(k-t+1)$, then,
$$
\tilde{M}(n,k,t)=M(n,k,t).
$$
\item
If $(t+1)(k-t+1)<n$ and $k\leq 2t+1$, then,
$$
\tilde{M}(n,k,t)=|\nu_1 (n,k,t)|.
$$
\item If $ (t+1)(k-t+1)<n$ and $k>2t+1$, then,
$$
\tilde{M}(n,k,t)=\max\{ |\nu_1 (n,k,t)|,|\nu_2 (n,k,t)|\} .
$$
\end{itemize}
\end{Th}
Note that the analogous value $M(n,k,t)$ was determined  in Ahlswede--Khachatrian \cite{11}. 
Before formulating our   main result, we make some additional definitions.
Set:
\begin{align}\label{e09}
 {\cal H}_i =\biggl\{ H\in{[t+i]\choose t+1}: [t]\subseteq H\biggr\} 
\bigcup ~\biggl\{ H\in{[t+i]\choose t+i-1}: [t+1,t+i]\subseteq H\biggr\}.
\end{align}
For ${\cal C}\subseteq 2^{[n]},$ we denote $\cal U({\cal C})$ as the set of permutations whose fixed point set is an upset of $\cal C$:
\[\cal U(\cal C)=\{p\in \Gamma: \exists\ \! C\in \cal C,~s.t.~ C\subseteq f(p) \}.\]
The main  result of this work is the following theorem. 
 \begin{Th}
 \label{th21}
 \begin{itemize}
 \item If, 
 $$
 \max\left\{ \ell=t+2r: \frac{\ell-t}{2(\ell -1)}\gamma (\ell ,n)\leq 1\right\} >t,
 $$
 then,
 $$
 \tilde{M}(n,t)=M(n,t).
 $$
 \item
 If,
 $$
 \max\left\{ \ell=t+2r: \frac{\ell-t}{2(\ell -1)}\gamma (\ell ,n)\leq 1\right\} =t ,
 $$
 then,
 $$
 \tilde{M}(n,t)=\max\{\nu_1 (n,t) ,\nu_2(n,t)\} ,
 $$
 where
 $$
 \nu_i (n,t)=\sum_{S\in {\cal U}({\cal H}_i )}\xi(n-|S|) .
 $$
 \end{itemize}
 \end{Th}
 
 Moreover, the result allows one to compute $\tilde M(n,t)$ for big $n$.
 
 \begin{cor}
 There is a constant $n_2(t)$ such that, for   $n>n_2 (t)$, 
 \[\tilde{M}(n,t)=(n-t)!-\xi (n-t)-\xi (n-t-1)+t. 
 \]
\end{cor}
 
 \section{$\cal F$- and $\cal L$-compressed sets}

We recall the fixing operation introduced in \cite{4}. For $i\neq j$ and $p\in\Gamma (n)$ we define the permutation 
$F(i,j,p)$ as
\begin{eqnarray}
\label{er09}
F(i,j,p)=\left\{\begin{array}{ll}
(p\setminus p_i )\bigcup \{ \{ i\} ,p(i) \setminus\{i\}\},& j = p(i),\\
p , \hbox{otherwise},
\end{array}
\right.
\end{eqnarray}
where $p_i =(i_1 ,\ldots ,i_m ,i , i_{m+1},\ldots ,i_{\ell})$ is the cycle from $p$ which contains $i$ and
$p_i \setminus i =(i_1 ,\ldots ,i_m ,i_{m+1},\ldots ,i_{\ell})$.
{We further define on ${\cal A}\subseteq\Omega (n,t)$}:
\begin{eqnarray}
\label{er08}
F(i,j,p,{\cal A})=\left\{\begin{array}{ll}
F(i,j,p),\ & F(i,j,p )\not\in {\cal A},\\
p, & F(i,j,p) \in{\cal A} .
\end{array}\right.
\end{eqnarray}
Finally, we define the \textit{fixing operation},
\begin{equation}
\label{er07}
{\cal F} (i,j,{\cal A})=\{ F(i,j,p,{\cal A}): p\in{\cal A} \}
\end{equation}
It is easy to see that ${\cal F}(i,j,{\cal A})$ preserves the cardinality of ${\cal A}$ and its $t$-intersection properties. Indeed, if a pair of permutations intersect in $t$ cycles which do not contain $i$,
then they will still intersect in these cycles after the application of fixing operations. 
Otherwise, if a pair of permutations intersect in a cycle that  contains $i$, the new permutations will intersect in the singleton $\{i\}$.
Moreover, the fixing operation preserves the volume of the family,  since \eqref{er09}-\eqref{er07} prevents that a permutation is obtained from two different permutations in the set. 

Repeatedly applying the fixing operation for different values of $i,j$ eventually results in a {set} ${\cal A}^*$ with the following property: for all $i,j\in [n],\ i\neq j$,  
$$
{\cal F}(i,j,{\cal A}^* ) ={\cal A}^*. 
$$
Moreover, as in \cite[Theorem 10]{4}, 
the fixed point sets of any pair of permutations $p_1 ,p_2 \in {\cal A}^*$ has at least $t$ common singletons. We state this fact as a lemma:
\begin{Le}\label{leF}
	Suppose that $\mathcal{A}$ is $\mathcal{F}$-compressed. Then, $\cal A$ is $t$-cycle intersecting if and only if
	\[ \Fix (\mathcal A)=\{f(p): p\in \mathcal A\}  \]
	is $t$-intersecting. Moreover, if there are $t$ cycles common to every permutation in $\cal A$, then there are $t$ singletons common to every element of $\cal A$.
\end{Le}
\begin{proof}
	Suppose there are permutations $p_1,p_2$ that intersect in $t$ cycles $\pi_1,...,\pi_t$. We claim that $\pi_1,...,\pi_t$ can be taken as singletons. Write $\pi_j=(i_1^j,...,i_{s_j}^j)$.  and $\cal F(i_1,i_2,p_2)$ are $t$-intersecting, 
	Since $\cal A$ is $\cal F$-compressed,  $p_1$ is $t$-intersecting with the permutation $p_2'$ resulting of successively applying every possible $\cal F(i_l^j,i_{l+1}^j,-)$ to $p_2$. Therefore, for each $\pi_i$, either $\pi_i$ is a singleton, or $p_1$ and $p_2$ shares an extra cycle  $\pi_{i}'$. In the last case, we can repeat the argument using $p_1$, $p_2'$ and $\pi_1',...,\pi_t'$, where $\pi_i'=\pi_i$ whether $\pi_i$ is  a singleton. On sees that a set of $t$ singletons, common to both $p_1$ and $p_2$,  is obtained by a finite number of the steps above.
	
	The last statement in the Lemma is proved by a similar argument,  where one replaces $p_1,p_2$ by any subset $p_1,...,p_r$ of permutations, and $p_1,p_2'$ by every possibility $p_1',p_2,...,p_r$;  $~~p_1,p_2',p_3,...,p_r$;...; $~~p_1,...,p_{r-1},p_r'$.
\end{proof}

In addition, if $\cal A$ is non-trivial $t$-cycle intersecting, $ \Fix (\cal A)$ is non-trivial $t$-intersecting. Moreover, if $\cal A$ is maximal in $\Omega(n,t)$ (resp. in $\tilde\Omega(n,t)$), $ \Fix (\cal A)$ completely defines $\cal A$:
\[ \cal A=\cal U( \Fix (\cal A)). \]
In particular, $ \Fix (\cal A)$ is maximal in its respective class: denote by $I'(n,t)$ the set of $t$-intersecting families of subsets of $[n]$ whose cardinality is not $n-1$. That is,
\[I'(n,t)=\{\cal A\in I(n,t): |A|\neq n-1,~\forall A\in \cal A  \}. \]
 Notice that $ \Fix (\cal A)\subseteq I'(n,t)$. One concludes that $ \Fix (\cal A)$ is maximal in $I'(n,t)$. In particular, $\cal F$ connects the problem of maximal families of permutations to the well-known problem of intersection of finite sets solved in \cite{7}. However, lifting from the fixed-point set to the family of permutations makes each step a little harder. The greatest contrast lies in the proof of Lemma \ref{l4}.

\bigskip 

Next, we define the (usual) {shifting procedure}.  Given $\ 1\leq v<w\leq n$ and a permutation with cycle decomposition
\[p=\{ ( j_1 ,\ldots ,j_{q-1},v,j_{q+1},\ldots ,j_s ),\ldots ,(w) ,\pi_1 ,\ldots ,\pi_c \}, \]
 define: 
$$
L(v,w,p)=\{ (j_1 ,\ldots ,j_{q-1},w,j_{q+1},\ldots ,j_s ) ,\ldots ,(v) ,\pi_1 ,\ldots ,\pi_c \} .
$$
Otherwise, if $p$ does not fix $w$, set,
$$
L(v,w, p)=p .
$$
We define $L(v,w, p, {\cal A})$ as follows:
$$
L(v,w,p,{\cal A})=\left\{\begin{array}{ll}
L(v,w,p),\ & L (v,w,p)\not\in{\cal A},\\
p, & L(v,w,p)\in{\cal A} .
\end{array}
\right.
$$
Finally the \textit{shifting operation} ${\cal L}(v,w,{\cal A})$ is defined as,
$$
{\cal L}(v,w,{\cal A})=\{ L(v,w, p,{\cal A}):\ p\in{\cal A}\}.
$$
It is clear that ${\cal L}(v,w,{\cal A})$ does not change the volume of ${\cal A}$, and preserves the $t$-cycle intersection property. Later we will show  that this operator preserves the nontrivial $t$-cycle intersection property as well (Lemma~\ref{s1}). Also, it is easy to see that after a finite number of operations we come  to a $\mathcal{L}$-\textit{compressed  $t$-intersecting set} ${\cal A}$, that is, 
$$
\mathcal L(v,w,{\cal A})={\cal A} ,\ 1\leq v<w\leq n.
$$
If $\mathcal A$ is $\mathcal F$-compressed as well, then each pair of permutations  has at least $t$ common singletons. In this case, the set of singletons $\Fix(\cal A)$ is \textit{left-compressed} (in the sense of \cite{7}). That is, if $\{a_1,...,a_k\}\in \Fix(\cal A)$, $a_1<a_2<...<a_j$, then $\{a_1',...,a'_k\}\in \Fix(\cal A)$ whenever $a_1'<a_2'<...<a_j'$, $a_i'\leq a_i$.

Next, we only consider the sets ${\cal A}$ that are both fixed by $\mathcal F$ and $\mathcal L$,
and denote the family of such  compressed  $t$-intersecting sets (respectively, non-trivial $t$-intersecting set) as $L\Omega (n,t)$ (respectively, $L\tilde{\Omega}(n,t))$. We begin the proof with the next lemmas.



Let $\Omega_0 (n,t)$ be the family of systems of permutations  ${\cal A}$, such that $|\bigcap_{p\in{\cal A}} f(p)|=0$.
\begin{Le}\label{s1}
	Let $\mathcal A\in \tilde\Omega(n,t)$ be such that $|\mathcal A|=\tilde M(n,t)$. Then, $\mathcal L(v,w,\mathcal A)\subset \tilde\Omega(n,t)$. In particular,
\begin{gather}
\label{e122}
\tilde{M}(n,t)=\max_{{\cal A}\in L\tilde{\Omega}(n,t)}|{\cal A}|.
\end{gather}
Moreover, if ${\cal A}\in\tilde{\Omega}(n,t)$ and $|{\cal A}|=\tilde{M}(n,t)$, then ${\cal A}\in\Omega_0 (n,t)$. In particular
\[M_0 (n,t)=\max_{{\cal A}\in\Omega_0 (n,t)}|{\cal A}|=\tilde{M}(n,t). \]
\end{Le}
\begin{proof}  First we will prove~(\ref{e122}). Let ${\cal A}\in\tilde{\Omega}(n,t)$, $|{\cal A}|=\tilde{M}(n,t)$, and suppose by contradiction that $L(v,w,{\cal A})\in\Omega (n,t)\setminus\tilde{\Omega}(n,t) .$ 
After  reordering $[n]$ and applying Lemma \ref{leF}, we can assume that $\cal A$ is $\cal F$-compressed, that $v=t, w=t+1$, and
\begin{gather}\label{eq1}
\bigcap_{p\in{\cal A}}f(p) =[t-1],\qquad  \bigcap_{p\in L(t,t+1,{\cal A})}f(p) =[t].
\end{gather}
In particular, if $[t]\nsubseteq f(p)$, $f(p)\cap[t+1]=[t-1]\cup \{t+1\}$. 
Therefore, $\cal A$ can be divided in three non-empty disjoint families: 
\begin{gather*}
\cal A_0=\{p\in \Gamma(n): [t+1]\subseteq f(p)\},\\
\cal A_1=\{p\in \cal A:  f(p)\cap [t+1]=[t]  \},\\
\cal A_2=\{p\in \cal A: f(p)\cap [t+1]=[t-1]\cup\{t+1\}  \}.
\end{gather*}

Observe that, if $ p\in \cal A_2$, then $ L(t,t+1,p)\notin \cal A$. Otherwise, $p\in L(t,t+1,\cal A)$ and, therefore,  the second equality in \eqref{eq1} is violated. On the other hand, using the maximality of $\cal A$, we claim that the simple permutations $p_1=(t+1,n)$, $p_2=(t,n)$ are in $\cal A$, reaching a contradiction with the aforementioned fact that $L(t,t+1,p_2)\neq \cal A$. From now on, we assume that $\cal L(v,w,\cal A)=\cal A$ for every  $t+1<v<w\leq n$. Observe that the last assumption preserves all properties described so far.

Suppose that $p_2\notin\cal A$. By observing that $|f(p)|>t$ for every $p\in \cal A$, we conclude that there is $p_n'$ such that $f(p_n')=[t]\cup \{n\}$. Since we assume that $\cal L(v,w,\cal A)=\cal A$ for every  $t+1<v<w\leq n$, there are permutations $p_{n-1}',~p_{n-2}',...,~p_{t+2}'$ satisfying  $f(p_i')=[t]\cup \{i\}$. Therefore, to preserve the $t$-intersection property between fixed point sets, we conclude that $f(p)\supset [t+2,n]$ for every $p\in \cal A_2$. But $|f(p)|=n-1$ in this case, violating the fact that permutations cannot move only one point. Following along the same lines, we conclude that $p_1\in \cal A$, reaching a contradiction with $L(t,t+1,\cal A)\notin\tilde\Omega(n,t)$.

We now prove the second part of the Lemma. Assume that ${\cal A}\in\tilde{\Omega}(n,t)\setminus\Omega_0 (n,t)$ and that $|{\cal A}|=\tilde{M}(n,t)$. We can also suppose that ${\cal A}\in L \tilde{\Omega}(n,t)$ and $1 \in f(p)$ for all $p\in{\cal A}$.
In particular, $p=(1,n) \notin \cal A$. Therefore, there exists $p_1 \in {\cal A}$ such that
$$
|[2,n-1]\bigcap f(p_1 )|\leq t-1.
$$
Following along the same lines as in the first part and using that $\cal A$ is $L$-shifted, we can conclude that there is $p_1\in \cal A$ such that $f(p_1 )=[t]\cup\{n\}$. In addition, all its shifts also appear: $f(p_2)=[t]\cup\{n-1 \}$, $f(p_3)=[t]\cup \{n-2\}$,..., $f(p_{n-t})=[t+1]$. We conclude that, if $p\in \cal A$, then either $f(p)\supset [t]$ or $f(p)\supset [t+1,n]$. Moreover, the latter case does not happen, since it must hold that $|f(p)\cap [t]|=t-1$, therefore $|f(p)|=n-1$, a contradiction. As a result, we conclude that  $[t]\subset f(p)$ for every $p\in \cal A$, contradicting  that $\cal A\in \tilde \Omega(n,t)$.
%
%
\end{proof}

\section{Generating sets and the proof of Theorem \ref{th21}}
Equipped with the preceding results, we are now ready to prove Theorem \ref{th21}. We follow along the lines of \cite{7} and denote $\Fix^* ({\cal A})$ as the set of minimal elements of $\Fix(\cal A)$, with respect to the set-theoretic inclusion. Recall that the inclusion ${\cal A}\in\Omega (n,t)\ (\cal A\in \tilde{\Omega}(n,t))$ is equivalent to the inclusion   $\Fix({\cal A})\in I(n,t)\ ( \Fix({\cal A})\in\tilde{I}(n,t))$. 
By using $\cal L$, we further assume that $\Fix({\cal A})$ is $\cal L$-compressed. Define,
\begin{eqnarray*}&&
s^{+}(\{a_1,....,a_j\}\subseteq[n])=\max_{i} a_i ,\\ &&
s^+ (\Fix({\cal A}))=\max_{a\in \Fix^* ({\cal A})} s^{+}(a) ,\\
&& s_{\min} =\min_{{\cal A}\in L\tilde{\Omega}(n,t): |{\cal A}|=\tilde{M}(n,t)}s^{+} (\Fix({\cal A})) .
\end{eqnarray*}
Given a maximal $\cal A$ and $C=f(p)$, $p\in \cal A$, we conclude that $p'\in \cal A$ whenever $f(p')\supseteq C$. Using that $\cal A$ is $\cal L$-shifted, we then conclude that
$$
{\cal A} =\bigcup_{C\in \Fix^* ({\cal A})} {\cal U}(\cal D(C)),
$$
where
$$
\cal D(C)=\left\{ A\in 2^{[n]} :\ A=C\cup B,  B\subseteq [s^{+} (C),n]\right\}
$$
We recall and adapt a couple of Lemmas from \cite{7}.
\begin{Le}
	\label{lem:AK4} Suppose that $\Fix^*(\cal A)$ is left-compressed and that $C\in g^* ({\cal A})$ is such that $ s^{+}(C)=s^{+}(\Fix({\cal A}))$.   Consider {the set of permutations generated by $C$ alone}
\begin{equation}
\label{e889}
{\cal A}_C =({\cal U}(C)\setminus {\cal U}(\Fix^* ({\cal A})\setminus\{ C\})).
\end{equation}
Then, $\Fix(\cal A_C)=\cal D(C)$.
\end{Le}

By $\cal A_C$, we mean that no other $E \in \Fix^* ({\cal A})$ is a subset of $\Fix(p)$ for any permutation  $p\in {\cal A}_C$.

\begin{Le}
	\label{lem:AK5} Consider $C_1 ,C_2 \in \Fix^* ({\cal A})$ where  $\cal A$ is left-compressed. If $i\not\in C_1 \cup C_2 ,\ j\in C_1 \cap C_2$ and $i<j$, then
	$$
	|C_1 \cap C_2 |\geq t+1.
	$$
\end{Le}
Lemmas \ref{lem:AK4} and \ref{lem:AK5} are straightforward  restatements of Lemmas 4 and 5 in \cite{7}, applied to singleton sets. 
The next Lemma helps us to establish the possible sets realized as $\Fix^* ({\cal A})$  for a maximal  ${\cal A}\in L\tilde{\Omega} (n,t)$
when $|{\cal A}|\neq M(n,t)$. 
The next Lemma is the main step in the proof of Theorem \ref{th21}. Observe that, although there are similarities, Lemma \ref{l4} does not follow from \cite[Lemma 6]{7}, since one cannot guarantee that $\Fix(\cal A)\nsubseteq {[n] \choose k}$ or that $|\Fix(\cal A)|=\tilde M(n,k,t)$, for any $k$. 
\begin{Le}
\label{l4}
Let  ${\cal A}\in L\tilde{\Omega} (n,t),\ |{\cal A}|=\tilde{M}(n,t)\neq M(n,t)$  be such
that $s^{+} (\Fix({\cal A}))=s_{\min} $. Then, for some $i\geq 2$, 
$$
\Fix^* ({\cal A})={\cal H}_i .
$$
\end{Le}
\begin{proof}
The proof follows closely the proof of \cite[Lemma 6]{7}, save the computation of $|\cal A|$ for specific $\cal A$'s. Let $\ell=s^{+}(\Fix({\cal A}))$ and consider
\[ \Fix_0 ({\cal A})=\{ g\in \Fix^* ({\cal A}):\ s^+ (g)=\ell\},\qquad \Fix_1({\cal A}) =\Fix^* ({\cal A})\setminus \Fix_0 ({\cal A}) .\]
It is easy to see that $\ell>t+1$.  It follows from Lemma \ref{lem:AK5} that, 
\begin{gather}\label{eq:P}
\text{if } C_1 ,C_2 \in \Fix_0 ({\cal A}),\ ~|C_1 \cap C_2 |=t \quad \Rightarrow \quad |C_1 |+|C_2 |=\ell +t . 
\end{gather}
Denote,
\[\biggl| \bigcap_{C\in \Fix_1 ({\cal A})}C\biggr| =\tau.\]

\begin{Cl}\label{Cl:1}
	$\tau\geq t$.
\end{Cl}
\begin{proof} Assume by contradiction that $\tau<t$. Consider the following partitioning of $\Fix_0(\cal A)$:
$$
\Fix_0 ({\cal A}) =\bigcup_{t<i<\ell} R_i ,\qquad R_i =\Fix_0 ({\cal A})\bigcap{[n]\choose i},
$$
and denote $$R^\prime_i =\{ C\subseteq [\ell -1]:\ C\cup \{\ell\}\in R_i \}.
$$
Since $\Fix({\cal A})$ is $\cal L$-compressed and $R_i'\subset {[n]\choose i}$, \eqref{eq:P} gives
  $$C_i \in R_i^\prime ,\ C_j \in R^\prime_j ,\ i+j\neq \ell +t\quad \Rightarrow \quad 
|C_i \cap C_j |\geq t.$$

Next, we show that $R_i =\emptyset.$ We claim that, if $R_j\neq \emptyset$ for some $j\neq (\ell+t)/2$, then there is a pair $R_i,R_{\ell+t-i}\neq\emptyset$, $i\neq (\ell+t)/2$. Assume, on the contrary, that for all $R_i \neq\emptyset$, $R_{\ell +t-i}=\emptyset $. But then, by \eqref{eq:P},
$$
\Fix^\prime =(\Fix^* ({\cal A})\setminus \Fix_0 ({\cal A}))\bigcup\bigcup_{t<i<\ell}R^\prime_i \in\tilde{I}(n,t)
$$
satisfy
$$
|{\cal U}(\Fix^\prime )|\geq |{\cal A}| \qquad \text{and} \qquad s^{+} (\Fix^\prime )<s^{+} (\Fix({\cal A})),
$$
which contradicts our initial assumption on $\cal A$. On the other hand, suppose there is  $i\neq (\ell+t)/2$ such that $R_i ,R_{\ell +t-i} \neq \emptyset .$ Consider the new sets
\begin{eqnarray*}
&&\varphi_1 =\Fix_1 ({\cal A})\bigcup \big(\Fix_0 ({\cal A})\setminus (R_i \cup R_{\ell +t-i})\big)\bigcup R^\prime_i ,\\
&&\varphi_2 =\Fix_1 ({\cal A})\bigcup \big(\Fix_0 ({\cal A})\setminus (R_i \cup R_{\ell +t-i})\big)\bigcup R^\prime_{\ell +t-i} .
\end{eqnarray*}
Since $\tau<t$, we have that $\varphi_i \in \tilde{I}(n,t),$ thus
$$
{\cal A}_i ={\cal U}(\varphi_i )\in\tilde{\Omega}(n,t) .
$$
This contradicts the maximality of $|\cal A|$ once we show that
\begin{equation}
\label{e56}
\max_{j=1,2}|{\cal A}_i |>|{\cal A}|.
\end{equation}
 Using~(\ref{e889}), we have:
\begin{eqnarray*}
&&|{\cal A}\setminus{\cal A}_1| =|R_{\ell +t-i}| \sum_{j=0}^{n-\ell}{n-\ell\choose j}\xi (n-\ell -t +i-j),\\
&&|{\cal A}_1 \setminus {\cal A}|\geq |R_i |\sum_{j=0}^{n-\ell}{n-\ell\choose j}\xi (n-i-j+1),\\
&& |{\cal A}\setminus{\cal A}_2 |=|R_i |\sum_{j=0}^{n-\ell}{n-\ell\choose j}\xi (n-i-j),\\
&&|{\cal A}_2 \setminus {\cal A}| \geq |R_{\ell +t-i}|\sum_{j=0}^{n-\ell}{n-\ell\choose j}\xi (n-\ell -t+i-j+1) .
\end{eqnarray*}
Therefore,  \eqref{e56} is violated only if both inequalities below are satisfied
\begin{align*}
|R_{\ell +t-i}| \sum_{j=0}^{n-\ell}{n-\ell\choose j}\xi (n-\ell -t +i-j)\geq |R_i |\sum_{j=0}^{n-\ell}{n-\ell\choose j}\xi (n-i-j+1),\\
|R_i |\sum_{j=0}^{n-\ell}{n-\ell\choose j}\xi (n-i-j)\geq |R_{\ell +t-i}|\sum_{j=0}^{n-\ell}{n-\ell\choose j}\xi (n-\ell -t+i-j+1).
\end{align*}
This is a contradiction, since $\xi(n+1)>\xi(n)$ for $n>0$. We thus conclude that $R_i =\emptyset $ when $\ i\neq (\ell +t)/2.$ 

Now assume $R_{\frac{\ell +t}{2}}\neq\emptyset.$ By the pigeon-hole principle, there exists a $k\in [\ell -1] $ and ${\cal S}\subseteq R^\prime_{(\ell +t)/2}$ such that $k\not\in B$ for every $ B\in{\cal S}$ and
\begin{equation}
\label{e9999}
|{\cal S}|\geq\frac{\ell -t}{2(\ell -1)}|R^\prime_{(\ell +t)/2}|.
\end{equation}
Hence, as before, we have $|B_1 \cap B_2 |\geq t$ for every $B_1 ,B_2 \in{\cal S}$, and
$$
D =(\Fix^* ({\cal A})\setminus R_{(\ell +t)/2})\bigcup {\cal S}\in\tilde{I}(n,t).
$$
Next we show that, 
\begin{equation}
\label{e4.9}
|{\cal U}(D )|> |{\cal A}|.
\end{equation}
Consider the partitions,
$$
{\cal A}={\cal G}_1 \cup{\cal G}_2,\qquad {\cal U}(D )={\cal G}_1 \cup{\cal G}_3,
$$
where
\begin{eqnarray*}
&&
{\cal G}_1 ={\cal U}(\Fix^* ({\cal A})\setminus R_{(\ell +t)/2}),\\
&&
{\cal G}_2 ={\cal U}(R_{(\ell +t)/2})\setminus {\cal U}(\Fix^* ({\cal A})\setminus R_{(\ell +t)/2})\\
&&
{\cal G}_3 ={\cal U}({\cal S})\setminus {\cal U}(\Fix^* ({\cal A})\setminus R_{(\ell +t)/2}).
\end{eqnarray*}
We will show that
\begin{equation}
\label{e4.10}
|{\cal G}_3 |>|{\cal G}_2 |.
\end{equation}
We have,
$$
|{\cal G}_2 |=|R_{(\ell +t)/2}|\sum_{j=0}^{n-\ell}{n-\ell\choose j}\xi \left( n-\frac{\ell +t}{2}-j\right),
$$
and,
$$
|{\cal G}_3 |\geq |{\cal S}|\sum_{j=0}^{n-\ell +1}{n-\ell +1 \choose j}\xi \left( n-\frac{\ell +t}{2}-j+1\right).
$$
Hence for~(\ref{e4.10}) to be true, it is sufficient that,
\begin{gather}\label{eq:proofgamma}
\frac{\ell -t}{2(\ell -1)}\gamma (\ell,n )>1.
\end{gather}We recall that, since we are assuming $\tilde M(n,t)\neq M(n,t)$, \eqref{eq:proofgamma} is always satisfied (see Theorem \ref{th00}. Observe that $\ell+t$ is even, otherwise $R_{\frac{\ell+t}{2}}=\emptyset$.) Hence $R_{\frac{\ell +t}{2}}=\emptyset .$
\end{proof}
Since $\cal A$ is $\cal L$-compressed, Claim \ref{Cl:1} gives
$$
\bigcap_{C\in \Fix_1 ({\cal A})} C=[\tau ]\qquad\text{and}\qquad
\ell =s^{+} (\Fix({\cal A}))>\tau. 
$$
Moreover, for all $C\in \Fix_0 ({\cal A})$,  $|C\cap [\tau ]|\geq\tau -1$ and, if $|C\cap [\tau ]|=\tau -1$, then $[\tau +1,\ell ]\subseteq C$.
Let us show that 
\begin{Cl}
	$\tau \leq t+1.$
\end{Cl}
\begin{proof}
	If $\tau\geq t+2$ then, for $C_1 , C_2 \in \Fix^* ({\cal A}),$
$$
|C_1 \cap C_2 \cap [\tau ] |\geq\tau -2\geq t.
$$
By denoting $\Fix^\prime_0 ({\cal A})=\{ C^\prime \subseteq [\ell -1]: C^\prime \bigcup\{\ell\}\in \Fix_0 ({\cal A})\}$, we have
$$
\varphi =(\Fix^* ({\cal A})\setminus \Fix_0 ({\cal A}))\bigcup \Fix^\prime_0 ({\cal A})\in\tilde{I}(n,t).
$$
Thus, 
 $|{\cal U}(\varphi )|\geq|{\cal A}|$ with  $s^{+}(\varphi )<\ell$, which contradicts the minimality of $\ell$.
\end{proof}

We are reduced to two cases:

If $\tau =t+1$ then $\ell =t+2$, otherwise, from the  argument above (by deleting $\{\ell \}$ from each element of $\Fix_0 ({\cal A})$) we obtain a generating set $\varphi \in\tilde{I}(n,t)$ satisfying $|{\cal U}(\varphi )|\geq |{\cal A}|$ and $s^{+}(\varphi )<\ell,$ a contradiction. On the other hand,  for these values of $\tau,\ell $ it clearly follows that  $\Fix^* ({\cal A})={\cal H}_2 .$

Finally, assume $\tau =t.$ Denote $\Fix^\prime_0 ({\cal A})=\{ C\in \Fix_0 ({\cal A}): |C\cap [t]|=t-1\} .$ 
We have
$$
\Fix^\prime_0 ({\cal A})\subseteq\{ C\subseteq [\ell ]:\ |C\cap [t]|=t-1, [t+1,\ell ]\subseteq C\} 
$$
and, for $C\in \Fix^* ({\cal A})\setminus \Fix_0^\prime ({\cal A})$, we have $[t]\subseteq C$ and $|C\cap [t+1,\ell ]|\geq 1.$ Hence, 
$$
{\cal U}(\Fix^* ({\cal A})) \subseteq {\cal U}({\cal H}_{\ell -t} ).
$$
Since ${\cal A}$ is maximal, we conclude that $\Fix^* ({\cal A})={\cal H}_{\ell -t} $, as desired.
\end{proof}

 Equipped with the preceding results, we are now ready to prove Theorem~\ref{th21}.
 
\begin{proof}[Proof of Theorem \ref{th21}]
	Denote $S_i =|{\cal U}({\cal H}_i )|$. It is only left to prove that the maximum value of $S_i$, with respct to $i$, is either achieved as $S_2$ or $S_{n-t-1}$. To this aim, we prove that, if
$S_i <S_{i+1}$, then $S_{i+1}<S_{i+2}.$
On the one hand, we have
$$
S_i =(n-t)!-\sum_{j=0}^{n-t-i}{n-t-i\choose j}\xi (n-t-j)+t\sum_{j=0}^{n-t-i}{n-t-i\choose j}\xi (n-t-i-j+1).
$$
On the other hand, $S_i<S_{i+1}$ if and only if
\begin{equation}
\label{e990}
\sum_{j=0}^{n-t-i-1}{n-t-i-1\choose j}\xi (n-t-j+1)\geq t\sum_{j=0}^{n-t-i-1}{n-t-i-1\choose j}\xi (n-t-j-i+1).
\end{equation}
Rewrite \eqref{e990} as follows:
\begin{gather*}
\sum_{j=0}^{n-t-i-2}{n-t-i-2\choose j}\xi (n-t-j+1)+\sum_{j=0}^{n-t-i-2}{n-t-i-2\choose j}\xi (n-t-j)\\
\geq t\sum_{j=0}^{n-t-i-2}{n-t-i-2\choose j}\xi (n-t-i-j+1)+t\sum_{j=0}^{n-t-i-2}{n-t-i-2\choose j}\xi (n-t-i-j) .
\end{gather*}
Therefore, it follows that:
\begin{equation}
\label{e91}
\sum_{j=0}^{n-t-i-2}{n-t-i-2\choose j}\xi (n-t-j+1)\geq t\sum_{j=0}^{n-t-i-2}{n-t-i-2\choose j}\xi (n-t-j-i),
\end{equation}
as desired. \end{proof}

As a last contribution, we fix $t$ and compute $\tilde{M}(n,t)$ for $n\to\infty$. First of all, observe that
$$
\frac{\sum_{j=0}^{n-t-2}{n-t-2\choose j}\xi\left(n-t-j\right)}{\sum_{j=0}^{n-t-2}{n-t-2\choose j}\xi\left( n-t-1-j\right)}=1+\gamma (t+2,n)\stackrel{n\to\infty}{\longrightarrow}\infty. 
$$
Moreover,
$$
\lim_{n\to \infty}\frac{\xi (n-t-1)}{\sum_{j=1}^{n-t-2}{n-t-2\choose j}\xi\left(n-t-j\right)}=0.
$$
It follows that,
for sufficiently large $n$,
\begin{align*}
S_{n-t}&= (n-t)!-\xi (n-t)-\xi (n-t-1)+t >S_{2}\\&=(n-t)!
-\sum_{j=0}^{n-t-2}{n-t-2\choose j}\xi (n-t-j)+t\sum_{j=0}^{n-t-2}{n-t-2\choose j}\xi (n-t-j-1),
\end{align*}
and hence, there is a constant $n_2(t)$ such that   $n>n_2 (t)$ implies  
$$
\tilde{M}(n,t)=(n-t)!-\xi (n-t)-\xi (n-t-1)+t .
$$

{\bf Aknowledgments}
\bigskip

The first author would like to express his deepest gratitude for his  collegues in Unifesp and USP, as well as the hospitality of both institutions.

\end{document}